\documentclass[letterpaper,10pt,journal,final,twocolumn]{IEEEtran}

\def\QED{\mbox{\rule[0pt]{1.3ex}{1.3ex}}}    
\def\endproof{\hspace*{\fill}~\QED\par\endtrivlist\unskip}

\usepackage{amssymb,amsmath,color}

\usepackage{enumitem}

\usepackage[textwidth=1.4cm,textsize=scriptsize]{todonotes}
\usepackage[english]{babel}
\usepackage{graphicx,epstopdf,float,tikz}
\usetikzlibrary{arrows,positioning,calc}

\usepackage[noadjust,compress]{cite}
\usepackage[normalem]{ulem}

\newtheorem{assumptionx}{Assumption}
\newtheorem{lemmax}{Lemma}
\newtheorem{propositionx}{Proposition} 
\newtheorem{remarkx}{Remark}
\newtheorem{theoremx}{Theorem}
\newtheorem{definitionx}{Definition}   
\newtheorem{examplex}{Example}

\newenvironment{assumption}{\smallskip\begin{assumptionx}\itshape}{\end{assumptionx}\smallskip}
\newenvironment{lemma}{\smallskip\begin{lemmax}\itshape}{\end{lemmax}\smallskip}

\newenvironment{theorem}{\smallskip\begin{theoremx}\itshape}{\end{theoremx}\smallskip}

\newcommand{\R}{\mathbb R}

\newcommand{\N}{\mathbb N}

\newcommand{\cA}{\mathcal A}

\newcommand{\cE}{\mathcal E}

\newcommand{\cN}{\mathcal N}

\newcommand{\cQ}{\mathcal Q}

\newcommand{\x}{\times}
\newcommand{\st}{\mid}

\newcommand{\und}{\underline}
\newcommand{\inv}{^{-1}}

\DeclareMathOperator{\linspan}{span}

\newcommand{\sr}{^\star}
 
\newcommand{\T}{^\top}
\newcommand{\eq}{^{\rm e}}
\newcommand{\m}{_{\rm m}}
\newcommand{\p}{_{\perp}}
\newcommand{\tPhi}{\widetilde{\Phi}}

\newcommand{\dx}{w}
\newcommand{\dz}{\nu}

\allowdisplaybreaks

\newcommand{\1}{\boldsymbol{1}}
\newcommand{\one}{\boldsymbol{1}}



\title{
\LARGE \bf Stability, Linear Convergence, and Robustness of the Wang-Elia  Algorithm  for Distributed Consensus Optimization
}
\author{Michelangelo Bin, Ivano Notarnicola and Thomas Parisini%  
}%%

\begin{document}

\onecolumn 
\vspace{4em}

\vspace{5em}
\begin{quote}
	\emph{\textcopyright{}~2022 IEEE.  Personal use of this material is permitted.  Permission from IEEE must be obtained for all other uses, in any current or future media, including reprinting/republishing this material for advertising or promotional purposes, creating new collective works, for resale or redistribution to servers or lists, or reuse of any copyrighted component of this work in other works.}
\end{quote}

\twocolumn
 
\maketitle

%%%%%%%%%%%%%%%%%%%%%%%%%%%%%%%%%%%%%%%%%%%%%%%%%%%%%%%%%%%%%%%%%%%%%%%%%%%%%%%%
\begin{abstract}
	We revisit an algorithm for distributed consensus optimization  proposed in 2010 by J.~Wang and N.~Elia. By means of a Lyapunov-based analysis, we prove input-to-state stability of the algorithm relative to a  closed invariant set composed of optimal equilibria and with respect to perturbations affecting the algorithm's dynamics. In the absence of perturbations, this result implies linear convergence of the local estimates and Lyapunov stability of the optimal steady state.  Moreover, we unveil fundamental connections with the well-known Gradient Tracking and with distributed integral control. Overall, our results suggest that a control theoretic approach can have a considerable impact on (distributed) optimization, especially when robustness is considered.
\end{abstract}

\section{Introduction}\label{sec.intro}

\subsection{Problem Overview and Literature Review}
We consider $N$ agents communicating through a connected, undirected  network represented by a simple graph $(\cN,\cE)$, with $\cN=\{1,\dots,N\}$ and $\cE\subset\cN^2$.
By exchanging information with neighbors, agents cooperatively seek a consensual solution $\theta\sr\in\R$ to the optimization problem\footnote{In this paper we focus on the single-variable case in which $\theta\in\R$. This simplifies the technical derivations without sacrificing generality, since all results reported in the paper directly extend  to the case where $\theta\in\R^m$, for some $m>1$, by properly introducing Kronecker products.} %Moreover, since this is anyway implied by     Assumption~\ref{ass.conv} introduced later, we assume from now on that $\theta\sr$ is unique.}
\begin{equation}\label{d.optimization_problem} 
\min_{\theta\in\R} \: \: \sum_{i\in\cN} f_i(\theta) 
\end{equation} 
where, for each $i\in\cN$, the function $f_i:\R\to\R$ is known to agent $i$ only.
Problem~\eqref{d.optimization_problem} is known as a \emph{cost-coupled} or \emph{consensus optimization} problem,  since agents minimize a global cost function $\sum_{i\in\cN}f_i$ over a common decision variable. As each agent $i\in\cN$ has only access to its own private function $f_i$, and not to the global cost function to be optimized, a distributed solution of Problem~\eqref{d.optimization_problem} is nontrivial.

Cost-coupled problems have been extensively investigated   in the last decades starting with the pioneering works~\cite{nedic2009distributed,nedic2010constrained,wang_control_2010}. 
A detailed account for the large amount of research on this topic can be found in the recent survey papers~\cite{molzahn2017survey,nedic2018network,nedic2018distributed,yang2019survey,notarstefano2019distributed}.
In particular,  an important step forward in the algorithmic solution  of \eqref{d.optimization_problem} was the introduction of a ``tracking'' protocol in the distributed gradient method. See, e.g.,~\cite{varagnolo2016newton,dilorenzo2016next,nedic2017achieving,
qu2018harnessing,xu2018convergence,
xi2018addopt,xin2018linear,
pu2020push,tian2020achieving} and the subsequent extensions~\cite{qu2019accelerated,xin2019distributed}. 
The algorithms based on this tracking protocol are known as \emph{Gradient Tracking} algorithms. According to the early interpretations, the tracking protocol aims at reconstructing, in a distributed way, the gradient of the global cost function. 
A   recent interpretation, instead, looks at the Gradient Tracking algorithms as embedding a distributed integral action~\cite{bin2019system}. As we discuss in Section~\ref{sec.main.GT}, this is one of the connection points with the Wang-Elia algorithm~\cite{wang_control_2010}  introduced later.

A main drawback of the Gradient Tracking algorithm is that it needs a specific initialization (see Section~\ref{sec.main.GT}) to work properly. As we clarify later in Section~\ref{sec.main.GT}, such initialization requirement makes the Gradient Tracking methods fragile with respect to uncertainties in the dynamics, such as those introduced by quantization, numerical errors in the computation of the gradients, or uncertainties affecting the communication with the neighbors. 
In particular, as the example in Section~\ref{sec.ex.GT} shows, even
a small quantization error can make the Gradient Tracking diverge to infinity, with a divergence rate that worsens for smaller stepsize values.

In~\cite{wang_control_2010}, a distributed algorithm was proposed for problem~\eqref{d.optimization_problem} that does not require any specific initialization. 
We refer to it as the \emph{Wang-Elia} algorithm. A continuous-time version of this algorithm was also studied in~\cite{hatanaka_passivity-based_2018} from a passivity-theoretic viewpoint. 
The discrete-time version, instead, represents the main subject of this work. 
In particular, in the Wang-Elia algorithm, each agent $i\in\cN$ maintains a pair of state variables $(x_i,z_i)\in\R^2$ that are updated as % follows
\begin{equation}\label{s.Wang-Elia.discrete}
\begin{aligned}
	x_i^+ &= x_i+ \!\sum_{j\in\cN_i} \!\beta a_{ij}\big(x_j-x_i + z_j-z_i\big)  - \alpha\beta\nabla f_i (x_i)\\
	z_i^+ &= z_i - \sum_{j\in\cN_i} \beta a_{ij}(x_j-x_i)   ,
\end{aligned}
\end{equation} 
in which $\cN_i:=\{j\in\cN\st (i,j)\in\cE\}$ is the neighborhood of $i$ (we stress that $i\notin \cN_i$) in the communication network $(\cN,\cE)$, $a_{ij}=a_{ji}>0$ for all  $(i,j)\in\cE$, and $\alpha,\beta>0$ are design parameters. The variable $x_i$ is the estimate Agent $i$ has of the optimal solution $\theta\sr$ of Problem~\eqref{d.optimization_problem},  and $z_i$ is an auxiliary state variable. 
%The result  of~\cite{wang_control_2010},  proved by means of input-output small-gain arguments (although the proof is only sketched)\todo{Not sure bout this}, states that, if each function $f_i$ is convex and $\beta,\alpha$ are chosen small enough, then all estimates $x_i$ converge to $\theta\sr$.
It was proved in~\cite{wang_control_2010} that, if each function $f_i$ is convex and  $\alpha,\beta$  are chosen small enough, then all estimates $x_i$ converge to $\theta\sr$.  

\subsection{Contribution}

We study the Wang-Elia algorithm~\eqref{s.Wang-Elia.discrete} in the presence of additive perturbations. 
We prove global \emph{input-to-state stability (ISS)}~\cite{sontag_smooth_1989} of the algorithm with respect to such perturbations and relative to a closed invariant set $\cA\sr\subset\R^{2N}$. The elements $(x,z)\in\cA\sr$ are all optimal in the sense that all estimates $x\sr_i$ equal the minimizer $\theta\sr$. In this way, we prove that, unlike the Gradient Tracking methods, the Wang-Elia algorithm is robust with respect to perturbations.
Moreover, in the absence of perturbation as in~\cite{wang_control_2010}, our results   establish Lyapunov stability of $\cA\sr$ and linear convergence of the local estimates $x_i$ to $\theta\sr$, which are stronger properties than only convergence as shown in~\cite{wang_control_2010}.
Finally, we compare the Wang-Elia and the Gradient Tracking algorithm, unveiling their similarities and differences, and making a connection with (distributed) integral control.
The developed analysis is based on Lyapunov arguments and provides further insights on the structure and functioning of the algorithm.
%
%Overall, our findings underline, once again, how a control-informed design/analysis approach may have a beneficial impact on (distributed) optimization and learning.

\subsection{Notation} \label{sec.intro.notation}
We denote 
% by $\R$ and $\N$ the real and natural numbers ($0\in\N$) and
by $\sigma(M)$ the spectrum of a matrix $M$ and we call it Schur if $\sigma(M)$ lies in the open unit disk. The vector and matrix-induced $2$-norms are denoted by $|\cdot|$. The distance of $x\in\R^n$, to a closed set $\cA\subset\R^n$ is denoted by $|x|_\cA:=\inf_{a\in \cA}|x-a|$. If $s:\N\to\R^n$,  we let $.^+$ denote the shift operator $s\mapsto s^+(\cdot)=s(\cdot+1)$, and $|s|_{t}:=\sup_{k=0,\dots,t}|s(k)|$. For compactness, we also  write $s^t$ in place of $s(t)$. For a given $N$, we let $\1:=(1,\dots,1)\in\R^N$ and we let $S\in\R^{N\x(N-1)}$ be a matrix satisfying 
\begin{align}\label{d.S}
S\T \one &=0, & S\T S &= I_{N-1}.
\end{align}
We define the matrix $T\in\R^{N\x N}$ and its inverse as
\begin{equation}\label{d.T}
T=\begin{bmatrix}
	\one\T/N \\S\T 
\end{bmatrix}, \qquad T\inv = \begin{bmatrix}
	\one& S
\end{bmatrix}.
\end{equation}
From \eqref{d.S}-\eqref{d.T}, we deduce that  the identity matrix $I_N$ satisfies
\begin{align}\label{e.I_one_S}
I_N &= {\one\one\T}/{N} + SS\T, & |S|&=1.
\end{align}
Given a $\chi\in\R^N$, we define its \emph{average-dispersion decomposition} as the pair $(\chi\m,\chi\p)=T\chi$, where $\chi\m:=\1\T\chi/N\in\R$  and $\chi\p:=S\T \chi\in\R^{N-1}$ are called, respectively, the \emph{average} and the \emph{dispersion} components of $\chi$. 
% From \eqref{d.S} and~\eqref{d.T}, we obtain 
Also, it holds $\chi = T\inv(\chi\m,\chi\p)= \1 \chi\m + S\chi\p$, and $|\chi|^2 = N\chi\m^2+|\chi\p|^2$.

\section{The Wang-Elia Algorithm Revisited}
\subsection{The Perturbed Wang-Elia  Algorithm}

In this paper, we study the following system 
\begin{equation}\label{s.xz.componentwise}
\begin{aligned}
	x_i^+ &= x_i \!+\! \sum_{j\in\cN_i} k_{ij}\big(x_j \!-\!x_i + z_j\!-\!z_i\big) \!-\! \gamma\nabla f_i (x_i) +\dx_i\\
	z_i^+ &= z_i \!-\! \sum_{j\in\cN_i} k_{ij}(x_j-x_i) +\dz_i  ,
\end{aligned}
\end{equation} 
for all $i\in\cN$, with arbitrary initial conditions $(x_i^0,z_i^0)$, with $\gamma$   a positive constant, and with $k_{ij}=k_{ji}>0$ for all $(i,j)\in\cE$. 
The terms $\dx:=(\dx_1,\dots,\dx_N)\!\in\R^N$ and $\dz:=(\dz_1,\dots,\dz_N)\!\in\R^N$ are perturbations modeling, e.g., uncertainties in the state measurements, in the computation of $\nabla f_i(x_i)$, and in the exchange of the neighboring states, or representing quantization errors and generic unmodeled dynamics.
The aggregate version of~\eqref{s.xz.componentwise} reads as% 
\begin{subequations}
\label{s.xz} 
\begin{align}
	x^+ &= (I-K)x - Kz -\gamma \Phi(x) + \dx, & x^0 & \in \R^N,
	\label{s.x}\\
	z^+ &= z + Kx + \dz, & z^0 & \in \R^N,
	\label{s.z}
\end{align} 
\end{subequations}
in which $K\in\R^{N\x N}$ is defined in such a way that  $K_{ij}=-k_{ij}$ for all $(i,j)\in\cE$, $K_{ii}= \sum_{j\in\cN_i} k_{ij}$ for all $i\in\cN$, and $K_{ij}=0$ otherwise, and
where $x:=(x_1,\dots,x_N)$, $z:=(z_1,\dots,z_N)$,  and $\Phi(x):=(\nabla f_1(x_1),\dots,\nabla f_N(x_N))$. 
Unlike \cite{wang_control_2010}, we do not factor $k_{ij}$ and $\gamma$ in terms of $\beta$ and $\alpha$ (cf.~\eqref{s.Wang-Elia.discrete}). We only assume that the coefficients $k_{ij}$ are chosen in such a way that $K$ satisfies the following conditions
\begin{align}\label{e.K}
K&=K\T, & \ker K &= \linspan\one, & \sigma(K)&\subset[0,1),
\end{align}
while the gain $\gamma$ is a small positive number to be chosen according to Theorem~\ref{thm.main} presented later in Section~\ref{sec.main.result}.

We underline that the last condition of~\eqref{e.K} is possible since the communication network is connected.  
Moreover,~\eqref{e.K} implies $\one\T K=0$ and that $S\T K S$ is invertible and Schur.

\subsection{Standing Assumptions}
We study System~\eqref{s.Wang-Elia.discrete} under the following assumptions. 
% Regarding the local functions $f_i$, we assume the following.
\begin{assumption}\label{ass.Lip}
For each $i\in\cN$, $f_i$ is continuously differentiable and $\nabla f_i$ is Lipschitz continuous.
\end{assumption} \vspace{-\smallskipamount}
%Finally, regarding Problem~\eqref{d.optimization_problem}, we assume the following.
\begin{assumption}\label{ass.conv}
The global cost function  $\sum_{i\in\cN} f_i$ is strongly convex.
\end{assumption}

Assumption~\ref{ass.conv} does not directly compare  to the assumptions of~\cite{wang_control_2010,hatanaka_passivity-based_2018}, where convexity of each $f_i$ is asked. Indeed, while Assumption~\ref{ass.conv} asks for strong convexity, such property is only required to the global cost function (as, e.g., in~\cite{tian2020achieving}), and not to each function $f_i$ individually. % as, e.g., in~\cite{wang_control_2010,qu2018harnessing,hatanaka_passivity-based_2018}.

Assumption~\ref{ass.conv} is not necessary to prove convergence of the estimates $x_i$ produced by~\eqref{s.xz}. However, when it holds, there is a natural choice among the optimal equilibria leading to a well-defined error system characterized by a simple structure. This supports a Lyapunov-based analysis allowing to establish, in addition to convergence, stronger stability and robustness properties.

\subsection{Existence of an Optimal Steady-State Locus}
Throughout the paper, when referring to an equilibrium of~\eqref{s.xz}, we always implicitly assume $(\dx,\dz)=0$.
We say that a state $(x,z)\in\R^{2N}$ is  \emph{consensually optimal} if $x_i=\bar\theta$ for all $i\in\cN$, where $\bar\theta\in\R$   is a critical point of the global cost function $\sum_{i\in\cN} f_i$ (i.e., $\sum_{i\in\cN}\nabla f_i(\bar\theta)=0$). The equilibria of~\eqref{s.xz} are characterized by the following lemma.
\begin{lemma}\label{lem.equilibria}
Suppose that $K$ in \eqref{s.xz} satisfies~\eqref{e.K}. Then, every equilibrium  of \eqref{s.xz} is consensually optimal. Conversely, if $\theta$ is a critical point of the global cost function, there exists an equilibrium $(x,z)$ of \eqref{s.xz} satisfying $x_i=\theta$ for all $i\in\cN$.
\end{lemma}

\begin{IEEEproof}
Consider \eqref{s.xz} with $(\dx,\dz)=0$. Then, $(x\eq,z\eq)$ is an equilibrium of \eqref{s.xz} if and only if
\begin{align} 
	Kx\eq+Kz\eq + \gamma\Phi(x\eq) &= 0, &
	Kx\eq &= 0.\label{pf.e.eq_conditions} 
\end{align} 
In view of~\eqref{e.K}, the second equation of~\eqref{pf.e.eq_conditions} is equivalent to $x\eq\in\linspan\1$. Hence,  $(x\eq,z\eq)$ is an equilibrium of \eqref{s.xz} if and only if $x\eq$ is a consensus point for the estimates $x_i$.

Let $\theta\eq\in\R$ be such that $x\eq=\1\theta\eq$. Then, the first claim follows by noticing that, in view of \eqref{e.K},  the first equation of~\eqref{pf.e.eq_conditions}  implies $\1\T\Phi(x\eq) = \textstyle\sum_{i\in\cN}\nabla f_i(\theta\eq) = 0$. 
%Namely, $\theta\eq$ is a stationary point of \eqref{d.optimization_problem} and, hence, $(x\eq,z\eq)$ is consensually optimal.

For the converse direction, let $\bar\theta\in\R$  be a stationary point of the global cost function, and let $x\eq:=\1\bar\theta$. Then, $\1\T\Phi(x\eq)=0$ and $Kx\eq=0$. Let $z\eq:= -\gamma S(S\T K S)\inv S\T \Phi(x\eq)$, where $S\T K S$ is invertible in view of \eqref{e.K}. Then, by repeatedly using \eqref{e.I_one_S}, and in view of~\eqref{e.K}, we get 
$Kz\eq 
= SS\T Kz\eq 
= -\gamma S (S\T K S)(S\T K S)\inv S\T \Phi(x\eq) 
= -\gamma SS\T \Phi(x\eq) 
= -\gamma  \Phi(x\eq)$. Hence, $(x\eq,z\eq)$ satisfies \eqref{pf.e.eq_conditions}.
\end{IEEEproof}
\smallskip

As the proof of Lemma~\ref{lem.equilibria} shows, the set of all equilibria of~\eqref{s.xz} (each of which is consensually optimal) can be expressed as follows
\begin{align*}
\cA\sr\!\!=\! \Big\{\! (x,z)\!\in\R^{2N}\!\!\st &\, \exists \theta\in\R,\, \textstyle\sum_{i\in\cN}\!\nabla f_i(\theta)=0,\, x=\1\theta,\\& \, z\!\in\!- \gamma S(S\T\! K S)\inv\! S\T \Phi(\1\theta) \!+\! \linspan\!\1  \Big\}, 
\end{align*}  
which is closed but not compact. We point out that Lemma~\ref{lem.equilibria} does not rely on the smoothness  and convexity assumptions.

The set $\cA\sr$ is the target steady-state locus of the forthcoming stability results and analysis. We stress that we cannot target a compact subset of $\cA\sr$ if global convergence is sought. Indeed, in view of \eqref{e.K}, even with $(\dx,\dz)=0$ the average component $z\m=\1\T z/N$ of $z$ remains constant along every solution of \eqref{s.xz}. We stress that the same holds also for the original algorithm~\eqref{s.Wang-Elia.discrete} as well as for the continuous-time counterpart, which therefore cannot have a compact attractor. We   underline that this property holds also for the Gradient Tracking algorithm, see Section~\ref{sec.main.GT}. 

\subsection{Main Result and Discussion}\label{sec.main.result}

A tuple $(x,z,\dx,\dz):\N\to \R^{4N}$ satisfying \eqref{s.xz} is called a \emph{solution tuple} of \eqref{s.xz}. We say that $\dz$ is \emph{integral-average bounded} if $t\mapsto \sum_{t\in\N} \dz\m(t)$ is bounded, where $\dz\m$ denotes the average component of $\dz$ (Section~\ref{sec.intro.notation}).

\begin{theorem}\label{thm.main}
Suppose that Assumptions~\ref{ass.Lip} and \ref{ass.conv} hold and that $K$ in \eqref{s.xz} satisfies \eqref{e.K}. Then, there exist $\gamma\sr,\alpha>0$  and, for each $\gamma\in(0,\gamma\sr)$, there exist $\mu_\gamma \in[0,1)$ and $\rho_\gamma,\tau_\gamma >0$, such that, for all $\gamma\in(0,\gamma\sr)$, every solution tuple $(x,z,\dx,\dz)$ of \eqref{s.xz} satisfies
\begin{equation}\label{e.claim.iss}
	\begin{aligned}
		|(x^t\!,z^t)|_{\cA\sr} &\!\le\! \alpha\mu_\gamma^t |(x^0\!,z^0)|_{\cA\sr}\!+\! \rho_\gamma |\dx\m|_{t-1} \!+\! \tau_\gamma|(\dx\p\!,\dz\p)|_{t-1}
	\end{aligned}
\end{equation}
for all $t\in\N$. In particular, if $\dx$ and $\dz$ are bounded, then  $x$ and $z\p$ are bounded. 
Moreover, if and only if $\dz$ is integral-average bounded, also $z\m$ (hence, $(x,z)$) is bounded.
\end{theorem}

Theorem~\ref{thm.main} is proved in Section~\ref{sec.analysis}. Under Assumption~\ref{ass.conv}, $|x-\1\theta\sr|\le|(x,z)|_{\cA\sr}$.
Hence, when $(\dx,\dz)=0$, Theorem~\ref{thm.main} implies exponential convergence of the estimates $x_i$ to the optimum~$\theta\sr$ with convergence rate $\mu_\gamma=\sqrt{1-c_0\gamma}$, being $c_0$ related to the convexity parameter of the global cost function (see Section~\ref{sec.analysis}).
We stress that convergence is \emph{global} in the initial conditions, unlike the Gradient Tracking (see Section~\ref{sec.main.GT} below).
Moreover, by means of standard ISS arguments~\cite{sontag_smooth_1989}, one can show that \eqref{e.claim.iss} implies
\begin{equation*}
\limsup |(x^t,z^t)|_{\cA\sr} \le \limsup\big( \rho_\gamma  |\dx\m^t| + \tau_\gamma   |(\dx\p^t,\dz\p^t)|\big).
\end{equation*} 
Thus, in particular, the estimates converge to $\theta\sr$ at front of every vanishing perturbation.

Furthermore, Theorem~\ref{thm.main} implies that the set $\cA\sr$ is Lyapunov stable when $(\dx,\dz)=0$, and strongly stable when  $(\dx,\dz)\ne0$. Namely, for every $\varepsilon>0$, there exists $\delta_{\varepsilon}>0$ such that $\max\big\{|(x^0,z^0)|_{\cA\sr},\, \sup_{t\in\N}|\dx\m^t|,\, \sup_{t\in\N}|(\dx\p^t,\dz\p^t)|\big\}<\delta_{\varepsilon}$ implies $|(x^t,z^t)|_{\cA\sr}<\varepsilon$ for all $t\in\N$. Nevertheless, we stress that the average component $z\m$ of $z$ may become unbounded when $\dz\m$ is  not integral-average bounded even if $|(x^t,z^t)|_{\cA\sr}\to 0$ and $\dz\m$ is small. Indeed, $z\m$ is Lyapunov stable when $\dz\m=0$ but not strongly stable when $\dz\m\ne 0$. It is, however, integral-ISS~\cite{angeli_characterization_2000} as established by Theorem~\ref{thm.main}.

Regarding the asymptotic gain property,  we underline that, as shown in Section~\ref{sec.analysis}, the gain $\rho_\gamma$ is $O(\gamma\inv)$ and $\tau_\gamma$ is $O(\gamma^{-1/2})$. Hence, the effect of $\dx\m$ and $(\dx\p,\dz\p)$ is, in general, amplified by taking smaller values of $\gamma$ ($\dz\m$, instead, is unaffected by $\gamma$). Nevertheless, in the relevant case where $\dz=0$ and $\dx$ represents uncertainty in the computation of the gradients $\Phi(x)$, we  have $w=\gamma w'$ for some $w'$, as this gives the term $-\gamma(\Phi(x)+w')$ in~\eqref{s.x}. In this case, the gain from $w'$ to $|(x^t,z^t)|_{\cA\sr}$ is $O(1)$.

Finally, we remark that the proof of Theorem~\ref{thm.main} is based on a \emph{time-scale separation}, enforced when $\gamma\ll 1$, between the average and the dispersion dynamics. In particular, the dynamics governing the consensus error is fast, while convergence of the average to the optimum is slow. In Section~\ref{sec.analysis}, these two dynamics are at first studied separately, and then   interconnected (see Figure~\ref{fig.interconnection}). It is interesting to notice that, while establishing stability of the average dynamics alone does put some constraints on $\gamma$, the condition $\gamma<1$ that  actually separates the time scales only arises when the two dynamics are interconnected.\footnote{In particular, $\gamma\le\gamma\sr$ in Section~\ref{sec.proof.interconnection} implies $\gamma\le(2c_7)\inv\implies\gamma\le1/(2\sqrt{N})<1$, being $c_7=c_2+c_4\ge c_4 =\sqrt{N}$.}

\section{Connections with the Gradient Tracking}\label{sec.main.GT}
%In its original formulation~\cite{varagnolo2016newton,dilorenzo2016next,nedic2017achieving,
%  qu2018harnessing,xu2018convergence,xi2018addopt,xin2018linear,pu2020push,tian2020achieving}, 
%	the Gradient Tracking algorithm employs a pair $(x_i,s_i)\in\R^{2N}$ of variables for each agent $i\in\cN$, whose (aggregate) update law reads as follows 
%\begin{equation}\label{s.GT_original}
%	\begin{aligned}
%		x^+ &= R x - \gamma s, & x^0&\in\R^N\\
%		s^+ &= Cs +\Phi(x^+) - \Phi(x), & s^0 &= \Phi(x^0),
%	\end{aligned}
%\end{equation} 
%where $R$ (resp.~$C$) is a row (resp.~column) stochastic matrix compatible with the communication network $(\cN,\cE)$.
%As for the Wang-Elia algorithm, convergence of the estimates $x_i$ to $\theta\sr$ is obtained under Assumptions~\ref{ass.connected}, \ref{ass.Lip}, and \ref{ass.conv} at an exponential rate. An equivalent (up to a change of coordinates) formulation of~\eqref{s.GT_original} was proposed in~\cite{bin2019system} and read as
%\begin{equation}\label{s.GT_canonical}
%	\begin{aligned}
%		x^+ &= R x + z - \gamma \Phi(x), & x^0 &\in \R^N\\
%		z^+ &= Cs  - \gamma(C-I)\Phi(x), & \1\T z^0 &= 0.
%	\end{aligned}
%\end{equation}
%Compared to~\eqref{s.GT_original}, the algorithm~\eqref{s.GT_canonical} is causal and does not need the gradients $\Phi(x^0)$ for the initialization.

In the ``canonical coordinates'' formulation\footnote{System~\eqref{s.GT_canonical} differs from the original formulation of the Gradient Tracking (see, e.g., ~\cite{varagnolo2016newton,dilorenzo2016next,nedic2017achieving,
	qu2018harnessing,xu2018convergence,xi2018addopt,xin2018linear,pu2020push,tian2020achieving}) by a change of coordinates and it is therefore equivalent. Nevertheless,~\eqref{s.GT_canonical} is causal and has the advantage of not requiring the computation of $\nabla f_i$ for the initialization.} of~\cite{bin2019system},
the Gradient Tracking algorithm employs a pair $(x_i,z_i)\in\R^{2}$ of variables for each agent $i\in\cN$, whose (aggregate) update law reads as follows:
\begin{subequations}\label{s.GT_canonical}
\begin{align} 
	x^+ &= R x + z - \gamma \Phi(x), && x^0 \in \R^N,\\
	z^+ &= Cz   - \gamma(C-I)\Phi(x), && \1\T z^0 = 0,\label{s.GT_canonical.z} 
\end{align} 
\end{subequations}
in which $R\in\R^{N\x N}$ (resp.~$C\in\R^{N\x N}$) is a row (resp.~column) stochastic matrix matching the communication network $(\cN,\cE)$, i.e., $R_{ij}=0$ (resp. $C_{ij}=0$) if $(i,j)\notin \cE$.
Like algorithm~\eqref{s.xz}, convergence to $\theta\sr$ of the estimates $x_i$ produced by~\eqref{s.GT_canonical} is obtained, at an exponential rate, under Assumptions~\ref{ass.Lip}, and~\ref{ass.conv}.

It is interesting to compare the Gradient Tracking~\eqref{s.GT_canonical} to algorithm~\eqref{s.xz} considered here.
First, we notice that also in~\eqref{s.x} the matrix $I-K$ multiplying $x$ is row stochastic in view of~\eqref{e.K}. Indeed, it is doubly stochastic. Likewise, the identity matrix multiplying $z$ in~\eqref{s.z} is column stochastic, and the exogenous term $Kx$ sums to zero as $- \gamma(C-I)\Phi(x)$ does in \eqref{s.GT_canonical}. Indeed, this implies that, like  algorithm~\eqref{s.xz}, also the Gradient Tracking has the property that $z\m=\1\T z/N$ is constant along every solution. Hence, the need of the initialization $\1\T z^0=0$ in~\eqref{s.GT_canonical}, which is the most significant difference between~\eqref{s.xz} and~\eqref{s.GT_canonical}.
As clear from the analysis in Section~\ref{sec.analysis} (see, in particular, Equations~\eqref{e.TKTinv} and~\eqref{s.error}), a similar initialization is not required for~\eqref{s.xz} because the uncontrolled dynamics $z\m$ is decoupled from the  other components of~\eqref{s.xz}.
We notice, indeed, that~\eqref{e.K} implies $Kz=KSz\p$. Hence,   $z\m$ is always filtered out in~\eqref{s.x}.

As for what concerns robustness, we underline that the unavoidable initialization and the coupling of $z\m$  with the remaining states make the Gradient Tracking~\eqref{s.GT_canonical} fragile if disturbances are added as in~\eqref{s.xz}.
Indeed, like in~\eqref{s.xz}, the uncontrolled dynamics $z\m$ of the Gradient Tracking can be destabilized by means of a bounded yet arbitrarily small additive perturbation $\dz$. 
However, unlike~\eqref{s.xz}, in the case of the Gradient Tracking $z\m$ affects all the other state variables. Hence, in general, an ISS result as that established by Theorem~\ref{thm.main} cannot not hold for~\eqref{s.GT_canonical}. A counterexample in this direction is given in Section~\ref{sec.ex.GT}.

Finally, we notice that, when $\dz=0$, Equation~\eqref{s.z} takes the form of an integrator processing the term $Kx$. From~\eqref{s.z}, by using~\eqref{e.K}, we can derive the following equation for the dispersion  component $z\p$ of $z$
\begin{equation}\label{s.zp}
z\p^+ =   z\p + S\T K  (x-\1 x\m).
\end{equation}

Since $S\T K$ is full row rank,~\eqref{s.zp} is an integral action processing the consensus error $x-\1 x\m$. Therefore, the Wang-Elia algorithm can be seen as a distributed proportional-integral (PI) controller (the proportional part being $(I-K)x-\gamma\Phi(x)$ and the integral part $Kz=KSz\p$) regulating the ``plant'' $x^+=u$ to the optimal equilibrium $\1\theta\sr$. 

Interestingly, it can be shown that the same distributed PI structure is shared also by the Gradient Tracking algorithm~\eqref{s.GT_canonical}, where the integrator processes the term $(C-I)(R-I)(x-\1 x\m)$ and only shows up in the coordinates $(x,z)\mapsto (x,(C-I)x-z)$.
However, it is worth noticing that, differently from the Gradient Tracking, the additional dynamics $z\m$ never contributes to the PI controller in~\eqref{s.xz}, regardless of how $z$ is initialized. Nevertheless, it still plays a crucial role since it enables the distributed implementation of the integral action otherwise impossible. In fact,~\eqref{s.zp} cannot be implemented in a distributed way since $S\T K$ does not match the sparsity constraints imposed by the communication structure.
%
%Hence, we conclude that both~\eqref{s.xz} and the Gradient Tracking can be simply seen as different, yet very similar, distributed PI controllers, the main difference being the gain matrices and the entailed initialization constraints. 

% \todo{OLD:}
% \MB{The additional dynamics $z\m$ does not contribute to the PI controller. Nevertheless, it plays a crucial role since it enables the distributed implementation of the integral action otherwise impossible. In fact,~\eqref{s.zp} cannot be implemented in a distributed way since $S\T K$ does not match the sparsity constraints imposed by the communication structure.
% 
% Interestingly, it can be shown that the same distributed PI structure is shared also by the Gradient Tracking algorithm~\eqref{s.GT_canonical} where, however, the integrator processes the term $(C-I)(R-I)(x-\1 x\m)$, and it only shows up in the coordinates $(x,\eta)$ where $\eta:=(C-I)x-z$.
%  Hence, we conclude that both~\eqref{s.xz} and the Gradient Tracking can be simply seen  as different, yet very similar, distributed PI controllers, the main difference being the gain matrices  and the entailed initialization constraints.}

\section{Stability Analysis}\label{sec.analysis}
In this section,  prove  Theorem~\ref{thm.main}. For ease of exposition, the proof is split in four parts.

\subsection{The  Reduced Error Subsystem}
Under Assumption~\ref{ass.conv}, there exists a unique $\theta\sr\in\R$ such that $(x,z)\in\cA\sr$ if and only if $x=\1\theta\sr$ and $z\in -\gamma S(S\T K S)\inv S\T \Phi(\1\theta\sr)  + \linspan\1$. Thus, we can define without ambiguity the equilibrium $(x\sr,z\sr)$ as
\begin{align*}
x\sr&:=\1\theta\sr, & z\sr&:=-\gamma S(S\T K S)\inv S\T \Phi(\1\theta\sr),
\end{align*}
and, with $T$ defined in \eqref{d.T}, change variables in~\eqref{s.xz} as
\begin{align*}
(x,z)\mapsto (\xi,\zeta) = (T(x-x\sr),T(z-z\sr)).
\end{align*}
These new variables represent the average-dispersion components (Section~\ref{sec.intro.notation}) of the errors $x-x\sr$ and $z-z\sr$. Indeed, $\xi=(\xi\m,\xi\p)$ and $\zeta=(\zeta\m,\zeta\p)$, with $\xi\m=\1\T(x-x\sr)/N=x\m-\theta\sr$, $\xi\p=S\T(x-x\sr)=S\T x=x\p$, $\zeta\m=\1\T(z-z\sr)/N = \1\T z/N = z\m$, and $\zeta\p=S\T(z-z\sr)$. In addition, we have
\begin{align*}%\label{e.xz_decompos}
x &= \1 (\xi\m+\theta\sr) + S\xi\p, & z&= \1 \zeta\m+S(\zeta\p + S\T z\sr ).
\end{align*}
The previous change of variables leads to the ``error system''  
\begin{align*}
\xi^{+} &= (I-T K T\inv) \xi  - T K T \inv \zeta    -\gamma T\widetilde\Phi(x)  + T\dx\\
\zeta^{+} &= \zeta + TKT\inv \xi + T\dz 
\end{align*} 
in which $\widetilde{\Phi}(x) :=\Phi(x)-\Phi(\1\theta\sr)$.
From   \eqref{d.T} and \eqref{e.K}, we get
\begin{equation}\label{e.TKTinv}
TKT\inv = \begin{bmatrix}
	\one\T K\one/N & \one\T K S/N\\
	S\T K\one & S\T KS
\end{bmatrix} \!= \!\begin{bmatrix}
	0 & 0 \\0 & S\T KS
\end{bmatrix}\!.
\end{equation} 
Hence,  the error system can be expanded as follows
\begin{subequations}\label{s.error}
\begin{align}
	\xi\m^{+} &= \xi\m -\gamma\1\T \tPhi(x)/N + \dx\m \label{s.error.xim}\\
	(\xi\p^{+},	\zeta\p^{+}) &= A (\xi\p, 
	\zeta\p)    -\gamma B \tPhi(x) + (\dx\p,\dz\p) 
	\label{s.error.xip_zetap}\\
	\zeta\m^{+} &= \zeta\m + \dz\m.\label{s.error.zetam}
\end{align} 
\end{subequations} 
in which
\begin{align}\label{d.AB}
A&:= \begin{bmatrix}
	I-S\T K S   & -S\T KS\\ 
	S\T KS &   I 
\end{bmatrix}, & B&:=\begin{bmatrix}
	S\T \\ 0
\end{bmatrix}.
\end{align} 

As clear from~\eqref{s.error},  the average component  $z\m=\zeta\m$  of $z$, which is marginally stable, is decoupled from the rest of the system. Indeed, $\zeta\m$ is not influenced by any other component of the state, nor it influences them. Moreover, under Assumption~\ref{ass.conv},  $\zeta\m$ does not contribute to the distance of $(x,z)$ to $\cA\sr$, as indeed we have 
\begin{equation}\label{e.dist_norm_relation}
\begin{aligned}
	|(x,z)|_{\cA\sr} 
	&= \inf_{(a,b)\in\cA\sr} |(x-a,z-b)| 
	\\
	&
	= \inf_{c\in\R} |(x-x\sr,z-z\sr+\1 c)|
	%		\\
	%	 	& =\inf_{c\in\R} \sqrt{N\xi\m^2 + |\xi\p|^2 + N(\zeta\m-c)^2+|\zeta\p|^2} 
	\\&= |(\sqrt{N}\xi\m,\xi\p,\zeta\p)|.
\end{aligned}
\end{equation}    
Therefore, we shall now drop Equation~\eqref{s.error.zetam} and focus on~\eqref{s.error.xim}-\eqref{s.error.xip_zetap}, to which we refer as the ``reduced error subsystem''. In the forthcoming Sections~\ref{sec.pf.xim} and~\ref{sec.pf.perp}, we analyze the two subsystems~\eqref{s.error.xim} and \eqref{s.error.xip_zetap} separately, and characterize their stability properties. Later  in Section~\ref{sec.proof.interconnection}, we   study their interconnection.

\subsection{The ``Average'' Subsystem $\xi\m$}\label{sec.pf.xim}
Define the function $V_1(\xi\m):=\xi\m^2$. The increment $\Delta V_1^t:=V_1(\xi\m^{t+1})-V_1(\xi\m^t)$ satisfies (here and in the following, we drop the time dependency when no confusion may arise)
\begin{subequations}\label{e.V1.increment} 
\begin{align}
	&\Delta V_1 = - 2\dfrac{\gamma}{N} \xi\m \one\T \tPhi(x)
	\label{e.V1.increment.1} 
	\\
	& \! +\!2\left( \xi\m \!-\!  \dfrac{\gamma}{N} \one\T\widetilde{\Phi}(x)\right)\! w\m 
	\!+\!|w\m|^2
	\!+\! \dfrac{\gamma^2}{N^2}|\1\T \tPhi(x)|^2 .  
	\label{e.V1.increment.2} 
\end{align} 
\end{subequations}
In view of strong convexity in Assumption~\ref{ass.conv}, we can write
\begin{equation}\label{e.convexity_stabilizes}
\begin{aligned}
	&\xi\m \1\T \tPhi(\1(\xi\m+\theta\sr)) 
	\\
	&= \xi\m \sum_{i\in\cN} \big(\nabla f_i(\xi\m+\theta\sr)-\nabla f_i(\theta\sr)\big) \ge 2c_0 N |\xi\m|^2
\end{aligned}
\end{equation}
for some $c_0>0$. Moreover, in view of Assumption~\ref{ass.Lip}, $\Phi$ is Lipschitz continuous and let $\ell$ be its Lipschitz constant. Then, by adding and subtracting $\gamma \xi\m\1\T\Phi(\1(\xi\m+\theta\sr))/N$ to~\eqref{e.V1.increment.1}, we obtain 
\begin{equation}\label{e.V1.t1.2}
\begin{aligned}
	\eqref{e.V1.increment.1} &=  - 2\dfrac{\gamma}{N} \xi\m \one\T\!\Big( \tPhi(\1(\xi\m\!+\!\theta\sr))\!  + \!\Phi(x)\! -\!\Phi(\1(\xi\m\!+\!\theta\sr)) \Big)\\
	&\le - 4c_0\gamma |\xi\m|^2 + 2 c_1 \gamma |\xi\m|\cdot |\xi\p|\\
	&\le -3c_0\gamma |\xi\m|^2 +  c_2 \gamma |\xi\p|^2,
\end{aligned}
\end{equation} 
with $c_1:=\ell /\sqrt{N}$, $c_2:=c_1^2/c_0$, and where  we used the Young's inequality
\begin{equation}\label{d.Youngs}
2ab \le \epsilon a^2 + b^2/\epsilon 
\end{equation}
with $a=|\xi\m|$, $b=|\xi\p|$, and $\epsilon= c_0/c_1$.
Similarly, we obtain
\begin{equation}\label{e.V1.t2.2}
\begin{aligned}
	& \eqref{e.V1.increment.2} \le   2\Big( (1+\gamma \ell)|\xi\m|+\sqrt{N} \xi\p|\Big) |w\m| + |\dx\m|^2\\
	&\qquad\qquad
	+ 
	2\ell^2\gamma^2 \left( |\xi\m|^2 +|\xi\p|^2/N\right) \\
	&\le \big(c_0\gamma/2 + c_3\gamma^2\big) |\xi\m|^2 + (c_4\gamma+c_5\gamma^2)|\xi\p|^2 + c_6(\gamma)|w\m|^2 
\end{aligned} 
\end{equation} 
in which $c_3:=c_0\ell/2 +2\ell^2$, $c_4 =\sqrt{N}$, 
$c_5:=2\ell^2/N$, $c_6(\gamma) = 1+2\ell/c_0 + (2/c_0+\sqrt{N})\gamma\inv$, and where we used~\eqref{d.Youngs} twice with $a=|\xi\m|$, $b=|w\m|$, $\epsilon=c_0\gamma/2$ and $a=|\xi\p|$, $b=|w\m|$, $\epsilon=\gamma$, respectively.
Pick $\gamma>0$ such that
\begin{equation*}
\gamma \le \gamma\sr_0 :=  c_0/(2c_3).
\end{equation*}
Then, \eqref{e.V1.t1.2} and \eqref{e.V1.t2.2} yield 
\begin{equation}\label{e.convexity_iss}
\Delta V_1 \le -2c_0\gamma |\xi\m|^2 + \big(c_7\gamma+c_5\gamma^2\big)|\xi\p|^2 + c_6(\gamma)|w\m|^2, 
\end{equation}
where $c_7:=c_2+c_4$.
The inequality~\eqref{e.convexity_iss} implies that the subsystem $\xi\m$ is exponentially ISS  relative to the origin and with respect to the inputs $\xi\p$ and $w\m$.

\subsection{The ``Dispersion'' Subsystem $(\xi\p,\zeta\p)$}\label{sec.pf.perp}
We now turn the attention to system~\eqref{s.error.xip_zetap}. First, we establish that the matrix $A$ in~\eqref{d.AB} is Schur, and hence that also~\eqref{s.error.xip_zetap} is ISS.
Let $\lambda\in\sigma(A)$ and $e=(e_1,e_2)$ a corresponding  unitary eigenvector. Then
\begin{equation}\label{e.A.lambda}
\lambda  = \lambda |e|^2 = e\T(\lambda e) = e\T A e 
= 1 - e_1\T S\T KS e_1.
\end{equation} 
Next, we claim that $e_1\ne 0$ for every eigenvector of $A$. Indeed, if $e=(0,e_2)$ and $\lambda\in\sigma(A)$, the equation $Ae=\lambda e$ implies $(S\T K S)e_2 = 0$, which implies $0\in\sigma(S\T K S)$ and thus contradicts \eqref{e.K}. Thus, since \eqref{e.K} also implies that $S\T K S$ is positive definite, we obtain from~\eqref{e.A.lambda} that $\lambda<1$ for all $\lambda\in\sigma(A)$. Finally,   $e_1\T S\T KS e_1\le \max\sigma(S\T KS) |e_1|^2< 1$, which together with \eqref{e.A.lambda} implies $\lambda>0$ for all $\lambda\in\sigma(A)$. Thus, $\sigma(A) \subset (0,1)$ and $A$ is Schur.

Let $\eta\p:=(\xi\p,\zeta\p)$ and $\delta\p:=(\dx\p,\dz\p)$. 
Define $V_2(\eta\p):=\eta\p\T P\eta\p$ with $P=P\T > 0$ being the unique solution to the Lyapunov equation $A\T PA-P = -3I$. The increment $\Delta V_2^t:=V_2(\eta\p^{t+1})-V_2(\eta\p^t)$ satisfies
\begin{subequations}
\begin{align}
	\Delta V_2 & = -3  |\eta\p|^2 -2\gamma \big(A\eta\p+\delta\p)\T  P B \tPhi(x)
	\label{e.V2.1}
	\\
	&
	+(2A\eta\p +\delta\p)\T  P\delta\p+ \gamma^2 \tPhi(x)\T B\T PB \tPhi(x).
	\label{e.V2.2}  
\end{align}
\end{subequations}

By using \eqref{d.Youngs} twice with $a=|\xi\m|$, $b=|\eta\p|$, 
$\epsilon=c_0/(4|A\T PB|\ell\sqrt{N})$ and $a=|\xi\m|$, $b=|\delta\p|$, $\epsilon=c_0/(4 |PB|\ell \sqrt{N})$ respectively, we obtain
\begin{align*}
\eqref{e.V2.1} 
\le 
\big( c_8 \gamma  - 3\big)|\eta\p|^2 + \dfrac{c_0 \gamma }{2} |\xi\m|^2 + c_9 \gamma |\delta\p|^2,
\end{align*}
with $c_8:=4|A\T PB|^2\ell^2 N /c_0 + 2|A\T PB|   \ell +|PB| \ell$ and 
$c_9:=4|PB|^2\ell^2N/c_0+|PB| \ell$.
Similarly, by using \eqref{d.Youngs} with $a=|\eta\p|$, $b=|\delta\p|$ and $\epsilon=1/(2|A\T P|)$, we obtain
\begin{equation*}
\eqref{e.V2.2} \le (1/2+c_{10} \gamma^2  )|\eta\p|^2 + c_{11}|\delta\p|^2 + c_{12}   \gamma^2 |\xi\m|^2 ,
\end{equation*}
where $c_{10}:=2|B\T P B| \ell^2$, $c_{11}:=2|A\T P|^2 + |P|$, and $c_{12}:=2|B\T P B| \ell^2 N$.
Pick $\gamma>0$ such that
\begin{equation*}
\gamma < \gamma\sr_1 := \min\big\{ \gamma\sr_0 ,\  (4 c_{10})^{-1/2},\ (4 c_8)\inv,\ c_0 (2 c_{12})\inv \big \}.
\end{equation*}
Then, with $c_{13}(\gamma):=c_{11} + c_9 \gamma $, it holds
\begin{equation}\label{e.DV2.iss}
\Delta V_2 \le -2|\eta\p|^2 + c_0 \gamma |\xi\m|^2 + c_{13}(\gamma) |\delta\p|^2.
\end{equation}
Similarly to~\eqref{e.convexity_iss},    inequality~\eqref{e.DV2.iss} establishes ISS of the dispersion subsystem  with respect to the average error $\xi\m$ and the dispersion component of the disturbances.

\subsection{The Interconnection Between $\xi\m$ and $(\xi\p,\zeta\p)$}\label{sec.proof.interconnection}

\begin{figure}[tpb]
\centering
\includegraphics[scale=.9,clip]{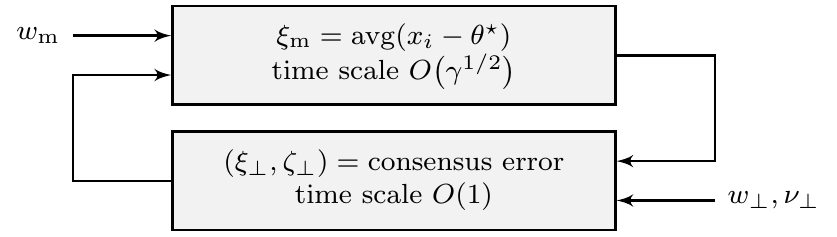}
% \vspace{-.15cm}
\caption{Feedback interconnection of $\xi\m$ and $(\xi\p,\zeta\p)$.}
\label{fig.interconnection}% \vspace{-.5cm}
\end{figure}

A block diagram representing the interconnection between~\eqref{s.error.xim} and~\eqref{s.error.xip_zetap} is represented in Figure~\ref{fig.interconnection} underlining the time-scale separation in the overall dynamics. 

Let $\und\lambda,\bar\lambda>0$ denote, respectively, the smallest and largest eigenvalues of $P$.
Pick $\gamma\in(0,\gamma\sr)$, where
\begin{equation*}
\gamma\sr := 
\min\big\{ \gamma\sr_1,\ (2c_7)\inv,\ (2c_5)^{-1/2},\ (c_0 \bar\lambda)\inv, \ c_0\inv \big\}.
\end{equation*}
% There follows from~\eqref{e.convexity_iss} and \eqref{e.DV2.iss} that, if $\gamma\in(0,\gamma\sr)$, where
Define the function $V(\xi\m,\eta\p) := V_1(\xi\m)+V_2(\eta\p) = |\xi\m|^2 + \eta\p\T P\eta\p$.
Then, in view of~\eqref{e.convexity_iss} and \eqref{e.DV2.iss}, the increment $\Delta V^t:=V(\xi\m^{t+1},\eta\p^{t+1}) - V(\xi\m^t,\eta\p^t)$ satisfies
\begin{equation}\label{e.V.iss}
\Delta V = - \gamma c_0 |\xi\m|^2 - |\eta\p|^2 + c_6(\gamma)|\dx\m|^2+ c_{13}(\gamma)|\delta\p|^2.
\end{equation} 
% where $\Delta V^t:=V(\xi\m^{t+1},\eta\p^{t+1}) - V(\xi\m^t,\eta\p^t)$. 
Equation~\eqref{e.V.iss} shows that the reduced error system is ISS with respect to the disturbances $\dx$ and $\dz\p$.

Then, notice that $|\xi\m|^2+\und\lambda |\eta\p|^2 \le V(\xi\m,\eta\p)\le |\xi\m|^2+\bar\lambda|\eta\p|^2$. Thus, \eqref{e.V.iss} implies
\begin{equation*}
V(\xi\m^t,\eta\p^t) \le q_\gamma^t V(\xi\m^0,\eta\p^0) + c_{14}(\gamma) |\dx\m|^2_{t-1} + c_{15}(\gamma)|\delta\p|^2_{t-1}
\end{equation*} 
where $q_\gamma := 1-\min\{c_0 \gamma ,\bar\lambda\inv\}  = 1-c_0 \gamma\in [0,1)$, $c_{14}(\gamma):=c_6(\gamma)/(c_0 \gamma)$ and $c_{15}(\gamma):=c_{13}(\gamma)/(c_0 \gamma)$. 

Since, in view of~\eqref{e.dist_norm_relation}, $(\max\{N,\und\lambda\inv\})\inv |(x,z)|_{\cA\sr}^2\le V(\xi\m,\eta\p)\le \max\{N\inv,\bar\lambda\}|(x,z)|_{\cA\sr}^2$, we finally obtain the sought inequality~\eqref{e.claim.iss} by setting 
\begin{align*}
\alpha & \! =\!  \sqrt{\max\{N\inv,\bar\lambda\} \max\{N,\und\lambda\inv\}},
\:\:
\mu_\gamma =\sqrt{q_\gamma},
\\
\rho_\gamma & \! =\! \sqrt{ c_{14}(\gamma) \max\{N,\und\lambda\inv\}},
% \\
\:\:
\tau_\gamma \! =\! \sqrt{c_{15}(\gamma) \max\{N,\und\lambda\inv\}}.
\end{align*}

Finally, the boundedness claims  directly follow    from~\eqref{e.claim.iss} and~\eqref{s.error.zetam}, respectively.\endproof

\section{Illustrative Example}\label{sec.ex.GT}
We present a toy example   showing the fragility of the Gradient Tracking.
We consider~\eqref{s.GT_canonical} for $N=2$ agents, with $C=R$, $R_{11}=R_{22}=0.8$, $R_{12}=R_{21}=0.2$ and $f_1(\theta):=(\theta-1)^2$, $f_2(\theta):=(\theta-4)^2$. We modify~\eqref{s.GT_canonical.z} to
%  with~\eqref{s.GT_canonical.z} substituted by
\begin{equation}\label{e.example.z}
z^+ = C \cQ (z) - \gamma(C-I)\Phi(x),
\end{equation}
where $\cQ (z):=10^{-5}\cdot \lfloor z\cdot10^{5} \rfloor$ models a quantization effect ($\lfloor \cdot \rfloor$ is the componentwise floor function). We can look at~\eqref{e.example.z} as the original~\eqref{s.GT_canonical.z} subject to the quantization error $\dz=C(\cQ (z)-z)$, which satisfies $\dz\m = \1\T\dz/N \le 0$.

\begin{figure}[tpb]
\centering
\includegraphics[scale=.9,clip]{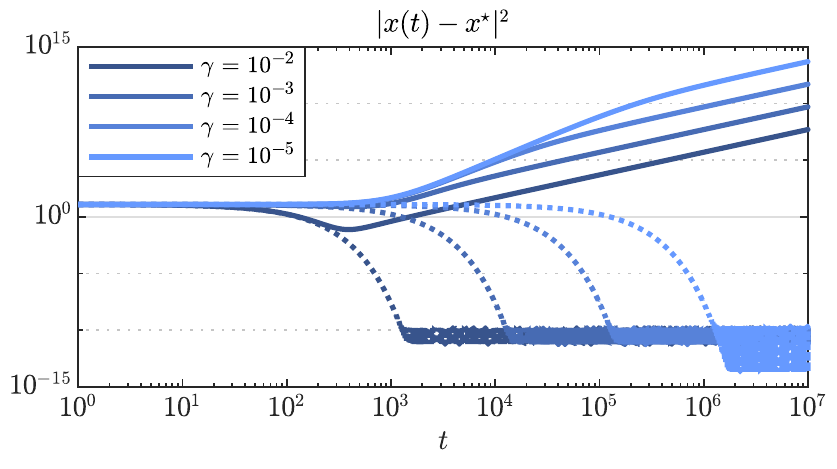}
%	\caption{Counterexample showing that the Gradient Tracking is not ISS.}
\vspace{-.2cm}
\caption{Comparison between algorithm~\eqref{s.GT_canonical} (solid lines) and algorithm~\eqref{s.xz} (dashed lines) showing the effect of the ISS property.}
\label{fig.GT}
% \vspace{-.4cm}
\end{figure}

Figure~\ref{fig.GT} shows four simulations obtained with stepsize $\gamma= 10^{-2} , 10^{-3} , 10^{-4}, 10^{-5}$ and with the same initial condition $(x^0,z^0)=(0,0)$. 
As discussed in Section~\ref{sec.main.GT}, the Gradient Tracking is not ISS. Indeed, the average quantization error $\dz\m$ destabilizes the state $x$. % , which escapes to the horizon. 
We stress that the smaller the $\gamma$ the higher is the divergence rate. 
This is explained by the same arguments given  in  Section~\ref{sec.main.result}.

Figure~\ref{fig.GT} also shows the solutions of algorithm~\eqref{s.xz} (for the same values of $\gamma$ and with the same initial condition), in which~\eqref{s.z} is modified to $z^+ = \cQ(z) + K x$.
Consistently with Theorem~\ref{thm.main}, we observe that the estimation error has a stable behavior, despite a small steady-state error.

\section{Conclusions}
We studied a perturbed version of the Wang-Elia algorithm, and we proved exponential ISS relative to an optimal steady state. 
%The analysis is based on Lyapunov arguments and highlights the presence of a time-scale separation between the consensus dynamics and the convergence of the average component to the optimum. 
We compared the algorithm to the Gradient Tracking, showing that the latter does not enjoy a similar ISS property due to the need of initialization. Overall, our arguments underline the impact that a control theoretic approach can have on the analysis of (distributed) optimization, especially when robustness is taken into account.

\bibliography{biblio,biblio_paper_GT}
\bibliographystyle{IEEEtran}

\end{document}